\documentclass[tran,10pt]{IEEEtran}

\usepackage{mathtools}
\usepackage[normalem]{ulem} 
\usepackage{subcaption}
\usepackage{enumerate}
\usepackage[ruled,vlined]{algorithm2e}
\usepackage{amsfonts}
\usepackage{amssymb}
\usepackage{booktabs}
\usepackage{multirow}
\usepackage{amsfonts,dsfont}

\usepackage[colorlinks=true,linkcolor=blue,citecolor=blue,urlcolor=blue,pdfborder={0 0 0},bookmarksnumbered]{hyperref}

\newcommand{\rev}[1]{{\color{blue}#1}}
\newcommand{\cut}[1]{}  

\newcommand{\vx}{{\mathbf{x}}}

\newcommand{\Diag}{{\mathrm{Diag}}} 

\DeclareMathOperator*{\argmin}{arg\,min}

\DeclareMathOperator*{\Min}{minimize}

\newcommand{\LPS}{\mathrm{LPS}}

\newcommand{\bc}{\begin{center}}
\newcommand{\ec}{\end{center}}

\newcommand{\beq}{\begin{equation}}
\newcommand{\eeq}{\end{equation}}

\newcommand{\bfl}{\begin{flushleft}}
\newcommand{\efl}{\end{flushleft}}

\newcommand{\bt}{\begin{tabbing}}
\newcommand{\et}{\end{tabbing}}

\newcommand{\beqn}{\begin{align}}
\newcommand{\eeqn}{\end{align}}

\newcommand{\beqs}{\begin{align*}} 
\newcommand{\eeqs}{\end{align*}}  

\newtheorem{theorem}{Theorem}

\newtheorem{definition}{Definition}

\makeatletter
\def\hlinewd#1{%
  \noalign{\ifnum0=`}\fi\hrule \@height #1 \futurelet
   \reserved@a\@xhline}

\usepackage{color,hyperref}
\definecolor{darkblue}{rgb}{0.0,0.0,0.4}
\hypersetup{colorlinks,breaklinks,
        linkcolor=darkblue,urlcolor=darkblue,
        anchorcolor=darkblue,citecolor=darkblue}

\makeatother
\title{Expander Graph and Communication-Efficient Decentralized Optimization}

\author{\IEEEauthorblockN{Yat-Tin Chow$^1$, Wei Shi$^2$, Tianyu Wu$^1$ and Wotao Yin$^1$\\}
        \vspace{0.4cm}
        \IEEEauthorblockA{$^1$Department of Mathematics, University of California, Los Angeles, CA\\
                $^2$ Department of Electrical \& Computer Engineering, Boston University, Boston, MA
}

\thanks{This work was supported in part by NSF grant ECCS-1462398.  Emails: \{ytchow,wuty11,wotaoyin\}@math.ucla.edu and wilburs@bu.edu}
        }
\renewcommand{\footnoterule}{%
        \kern -1pt
        \hrule width 1.1in height 1.2pt
        \kern 2pt
}

\IEEEoverridecommandlockouts
\begin{document}

\def\helvetica{phvr7t.tfm}
\def\helveticaoblique{phvro7t.tfm}
\def\helveticabold{phvb7t.tfm}
\def\helveticaboldoblique{phvbo7t.tfm}

\def\textcolor{ \footnote}

\def\rev{ }

=\helveticabold
=\helveticaboldoblique
\maketitle

\begin{abstract}
In this paper, we discuss how to design the graph topology to reduce the  communication complexity of certain algorithms for decentralized optimization. Our goal is to minimize the total communication needed to achieve a prescribed accuracy. We discover that the so-called \emph{expander graphs} are near-optimal choices. 

We propose three approaches to construct expander graphs for different numbers of nodes and node degrees. Our numerical results show that the performance of decentralized optimization is significantly better on expander graphs than  other regular graphs.

\end{abstract}

\section{Introduction}
\subsection{Problem and background}
Consider the decentralized consensus problem:
\begin{equation}
\Min_{x \in \mathbb{R}^P } f(x) = \frac{1}{M}
 \sum_{i= 1}^M  f_{i}(x)  \,,
\label{important}
\end{equation}
which is defined on a connected network of $M$ agents who cooperatively solve the problem for a common minimizer $x \in \mathbb{R}^P$. Each agent $i$ keeps its private convex function $f_i :  \mathbb{R}^P \rightarrow  \mathbb{R}$.
We consider synchronized optimization methods, where all the agents carry out their iterations at the same set of times.

Decentralized algorithms are applied to solve problem \eqref{important} when the data are collected or stored in a distributed network and a fusion center is either infeasible or uneconomical. 
Therefore, the agents in the network must perform local computing and communicate with one another. Applications of decentralized computing are widely found in sensor network information processing, multiple-agent control,  coordination, distributed machine learning, among others. Some important works include decentralized averaging \cite{EXTRA7, EXTRA15, EXTRA34}, learning \cite{EXTRA9, EXTRA22, EXTRA26}, estimation \cite{EXTRA1, EXTRA2, EXTRA16, EXTRA18}, sparse optimization \cite{EXTRA19, EXTRA35}, and low-rank matrix completion \cite{EXTRA20} problems.

In this paper, we focus on the communication efficiency of certain decentralized algorithms.  Communication is often the  bottleneck of distributed algorithms since it can cost more time or energy than computation, so reducing communication  is a primary concern in distributed computing. 

Often, there is a trade-off between \emph{communication requirement} and \emph{convergence rates} of a decentralized algorithm. Directly connecting more pairs of agents, for example, generates more communication at each iteration but also tends to make faster progress toward the solution. 
To reduce the \emph{total} communication for reaching a solution accuracy, therefore, we must find the balance. In this work, we argue that this balance is approximately reached by the so-called Ramanujan graphs, which is a type of expander graph.

\subsection{Problem reformulation}

We first reformulate problem \eqref{important} to include the graph topological information.  Consider an undirected graph $G = (V,E)$ that represents the network topology, where $V$ is the set of vertices and $E$ is the set of edges. Take a symmetric and doubly-stochastic \emph{mixing matrix} $W = \{ W_{e} \}_{e \in E}$, i.e. $ 0 \leq W_{ij} \leq 1$, $\sum_i W_{ij} = 1$ and $W^T = W$.  The matrix $W$ encodes the graph topology since $W_{ij} = 0$ means no  edge between the $i$th and the $j$th agents. 
Following \cite{shi2015extra}, the optimization problem \eqref{important} can be rewritten as:
\begin{equation}
\Min_{\textbf{x} \in \mathbb{R}^{M\times P} } \frac{1}{M} \sum_{i= 1}^M   F_{i}( \textbf{x} )~  \text{ subject to } W \textbf{x} = \textbf{x}  \,,
\label{min3}
\end{equation}
where $\vx=[x_1,...,x_M]^\top$ and the functions $F_{i}$ are defined as $ F_{i}( \textbf{x} ) := f_i(x_i) $.
\subsection{Overview of our approach}\label{sc:overview}

Our goal is to design a network of $N$ nodes and a mixing matrix $W$ such that state-of-the-art consensus optimization algorithms that apply multiplication of $W$ in the network, such as EXTRA \cite{shi2015extra,shi2015proximal} and ADMM \cite{ghadimi2014admm,shi2014linear,mota2013d}, can solve problem \eqref{important} with a nearly minimal amount of total communication. Our target accuracy is  $\|z^k-z^*\|<\varepsilon$, where $k$ is the iteration index, $\varepsilon$ is a threshold, and $(z^k)_k$ is the sequence of iterates, and  $z^*$ is its limit. We use $z^k,z^*$ instead of $x^k,x^*$ since EXTRA and ADMM iterate both primal and dual variables. \textbf{Let $U$ be the amount of \emph{communication per iteration}, and $K$ be the number of iterations needed to reach the specified accuracy.} The \emph{total communication} is their product $(UK)$. Apparently, $U$ and $K$ both depend on the network. A dense network usually implies a large $U$ and a small $K$. On the contrary, a highly sparse network typically leads to a small $U$ and a large $K$. (There are special cases such as the \emph{star} network.) The optimal density is typically neither  of them. Hence we shall find the optimal balance. To achieve this, we express $U$ and $K$ in the network parameters that we can tune, along with other parameters of the problem and the algorithm that affect total communication but we do not tune. Then, by finding the optimal values of the tuning parameters, we obtain the optimal network.

The dependence of $U$ and $K$ on the network can be complicated. Therefore, we make a series of simplifications to the problem and the network so that it becomes sufficiently simple to find the minimizer to $(UK)$. We mainly (but are not bound to) consider strongly convex objective functions and decentralized algorithms using fixed step sizes (e.g., EXTRA and ADMM) because they have the linear convergence rate:  $\|z^k-z^*\|\le C^k$ for $C<1$. Given a target accuracy $\varepsilon$, from $C^K=\varepsilon$, we deduce  the sufficient number of iterations $K=\log(\varepsilon^{-1})/\log(C^{-1})$. Since the constant  $C<1$ depends on network parameters, so does $K$. 

We mainly work with three network parameters: the condition number of the graph Laplacian, the first non-trivial eigenvalue of the incidence matrix, and the degree of each node. In this work, we restrict ourselves to the $d$-regular graph, e.g., every node is connected to exactly $d$ other nodes.

We  further (but are not bound to) simplify the problem by assuming that all edges have the same cost of communication. Since most  algorithms perform a fixed number (typically, one or two) of rounds of communication in each iteration along each edge,   $U$ equals a constant times $d$.

So far, we have simplified our problem to minimizing $(Kd)$, which reduces to  choosing $d$, deciding the topology of the $d$-regular graph, as well setting the mixing matrix  $W$. 

For the two algorithms EXTRA\ and ADMM, we deduce (in Section \ref{sc:commcost} from their existing analysis) that $K$ is determined by the mixing matrix $W$ and $\tilde{\kappa}(L_G)$, where $L_G$ denotes the graph Laplacian and $\tilde{\kappa}(L_G)$ is its reduced condition number. 
We simplify the formula of $K$ by relating $W$ to $L_G$ (thus eliminating $W$ from the computation).  
Hence, our problem becomes minimizing $(Kd)$, where $K$ only depends on $\tilde{\kappa}(L_G)$ and   $\tilde{\kappa}(L_G)$ further depends on the degree $d$ (same for every node) and the topology of the $d$-regular graph. 

By classical graph theories, a $d$-regular graph satisfies $\tilde{\kappa}(L_G) \le \frac{d+\lambda_1(A)}{d-\lambda_1(A)} $, where $\lambda_1(A)$ is the first non-trivial eigenvalue of the incidence matrix $A$, which we define later. We  minimize $\lambda_1(A)$  to reduce $\frac{d+\lambda_1(A)}{d-\lambda_1(A)}$ and thus $\tilde{\kappa}(L_G)$. Interestingly, $\lambda_1(A)$  is lower bounded by the Alon Boppana formula of $d$ (thus also depends on $d$). The lower bound is nearly attained by Ramanujan graphs, a type of graph with a small number of edges but having high connectivity and good spectral properties! Therefore, given a degree $d$, we shall choose Ramanujan graphs, because they approximately  minimize $\lambda_1(A)$, thus $\tilde{\kappa}(L_G)$, and in turn $(Kd)$ and $(KU)$. 

Unfortunately, Ramanujan graphs are  not known for all values of $(N,d)$. For large $N$, we follow \cite{friedman2008proof} to construct a random network, which has nearly minimal $\lambda_1(A)$.  For small $N$, we apply the Ramanujan Sparsifier \cite{batson2012twice}.
These approaches lead to near-optimal network topologies for a given $d$.

After the above deductions, minimizing the total communication simplifies to choosing the degree $d$. To do this, one must know the cost ratio between computation and communication, which depends on applications. Once these parameters are known, a simple one-dimensional search method will find $d$.

\subsection{Communication reduction}  \label{sec:reduction}

We find that, in different problems, the graph topology affects total communication to different extents. When one function $f_i(x)$ has a big influence to the consensus solution, 
 total communication is sensitive to the graph topology. On the other hand, computing the consensus average of a set of similar numbers is insensitive to the graph topology, so  a Ramanujan graph does not make much improvement. Therefore, our proposed graphs work only for the former type of problem and data.  
See section \ref{sec:num} for more details and comparison.
\subsection{Limitations and possible resolutions} \label{sec:limitations}
We assume nodes having the same degree,  edges having the same communication cost, and a special mixing matrix $W$. 
Despite these assumptions, the Ramanujan graph may still serve as the starting point to  reduce the total communication in more general settings, which is our future work.  

\subsection{Related work}

Other work that has considered and constructed Ramanujan graphs include \cite{kim2006maximizing,ghosh2006growing}.  
Graph design for optimal convergence rate under decentralized consensus averaging and optimzation is considered in \cite{EXTRA34,xiao2004fast,zhao2012diffusion,boyd2004fastest,yuan2013convergence,tsitsiklis1984problems,chen2012fast}.
\rev{However, we suspect the mixing matrices $W$ recommended in the above for average consensus might sometimes not be 
the best choices for decentralized optimization, since consensus averaging is a lot less communication demanding than decentralized optimization problems (c.f. Section \ref{sec:reduction}), and it is the least sensitive to adding, removing, or changing one node.
On the other hand, we agree for consensus averaging, those matrices recommended serve as the best candidates.} 

\section{Communication cost in decentralized optimization}\label{sc:commcost}

As described in Subsection \ref{sc:overview}, we focus on a special class of decentralized algorithm that allows us to optimize its communication cost by designing a suitable graph topology.

For the notational sake, let us recall several definitions and properties of a graph. See \cite{murty2003ramanujan,newman2000laplacian,chung1997spectral} for background. Given a graph $G$, the adjacency matrix $A$ of $G$ is a $|V| \times |V|$ matrix with components
\begin{equation*}
A_{ij} = \begin{cases}  1 \quad \text{ if } (i,j)\in E \,,  \\ 0 \quad \text{ if otherwise.} \end{cases}
\end{equation*}
An equally important (oriented) incidence matrix is defined as is a $|E| \times |V| $ matrix such that
\begin{equation*}
B_{e i} = \begin{cases}  1 \quad &\text{ if } e = (i,j) \text{ for some }j, \text{ where } j>i\,,  \\ -1 \quad &\text{ if } e = (i,j) \text{ for some }j,  \text{ where } j<i\,,  \\  0 \quad &\text{ otherwise.} \end{cases}
\end{equation*}
The degree of a node $i \in E$, denoted as $\deg(i)$, is the number of edges incident to the node.
Every node in a $d$-regular graph has the same degree $d$.
The graph Laplacian of $G$ is:
\begin{equation*}
L_G := D - A = B^T B \,,
\end{equation*}
where $D$ has diagonal entries $D_{ii} = \deg(i)$ and 0 elsewhere. 
The matrix $L_G$ is symmetric and positive semidefinite. The multiplicity of its trivial eigenvalue $0$ is the number of connected components of $G$.  Therefore, we consider its first \emph{nonzero} eigenvalue, denoted as $\widetilde{ \lambda_{\text{min}}} (L_G)$ and its largest eigenvalue, denoted as $\lambda_{\text{max}}(L_G)$.
The reduced condition number of  $L_G$ is defined as:
\begin{equation}
\tilde{\kappa}(L_G) := \lambda_{\text{max}}(L_G) / \widetilde{ \lambda_{\text{min}}} (L_G) \,,
\label{cond}
\end{equation} 
which is of crucial importance to our  exposition.
We  need a weighted graph and its graph Laplacian.  A weighted graph $G = (V,E,w)$ is a graph with nonnegative edge weights $w = (w_e)$.  Its graph Laplacian is defined as
\begin{equation*}
L_{G,w} := B^T \Diag(w) B \,.
\end{equation*}
Its  eigenvalues and condition numbers are similarly defined, e.g., $\tilde{\kappa}(L_{G,w})$ being its reduced condition number.

Next we define a special class of decentralized algorithm, which simplifies the computation of total communication cost.

\begin{definition}
\label{def:1}
We call a first-order decentralized algorithm to solve \eqref{min3} communication-optimizable if any sequence $\{ \textbf{z}^{k} \}$ generated by the algorithm satisfies
\begin{equation}
\label{R}
\|\textbf{z}^{k} - \textbf{z}^* \| \leq R( \Theta, \Gamma_f, \tilde{\kappa}(L_{G,w}), k  ) \,,
\end{equation}
for a certain norm $ \| \cdot \|$ and function $R$, where  $\Theta$ consists of all algorithmic parameters specified by the user, $\Gamma_f$ is a constant depending on the objective function $f$.
 In addition, $R$ is non-increasing as $\tilde{\kappa}(L_{G,w}) $ decreases or $k$ increases.
\end{definition}

The above definition means that there exists a convergence-speed bound that depends on only $L_{G,w}$ and $k$ while other parameters $\Theta, \mu_f$ are fixed.
This description of convergence rate is met by many algorithms mainly because $\tilde{\kappa}(L_{G,w})$ reflects the connectivity of graph topology.
In fact, if $G$ is $d$-regular and weighted, then we have from \cite{newman2000laplacian,chung1997spectral}:
\begin{equation*}
 \tilde{\kappa}(L_G)  \leq \frac{d + \lambda_1 (A) }{ d - \lambda_1 (A) },
\end{equation*}
where $\lambda_1 (A) \neq d $ is the first non-trivial eigenvalue of the adjacency matrix $A$ (as the eigenvalue $\lambda_0 (A)=d$ is trivial and corresponds to the eigenvalue $0$ of $L_G$) and is  related to the connectivity of $G$ by the Cheeger inequalities  \cite{alon2004probabilistic}:
\begin{equation*}
\tfrac{1}{2} ( d - \lambda_1 (A) ) \leq h(G) \leq \sqrt{ 2 d (d - \lambda_1 (A) )}
\label{eqneqneqn}
\end{equation*}
where $h(G)$ is the Cheeger constant quantifying the connectivity of $G$.  Therefore saying that $R( \Theta, \mu_f, \tilde{\kappa}(L_{G,w}), k  ) $ is nondecreasing with $\tilde{\kappa}(L_{G,w})$ is related to saying that is nonincreasing with the connectivity of $G$. (Higher connectivity lets the algorithm converge faster.) Hence, we  hope to increase the number of edges $E$ of $G$ though a cost is associated. 

With the help of  $L_{G,w}$, for simplicity we might choose the mixing matrix $W$ to take the form \cite{xiao2004fast,zhao2012diffusion}:
\begin{equation}
W_{L_{G,w}} := I - \frac{2}{ (1+ \Theta_1) \lambda_{\text{max}}(L) } L_{G,w} \, ,  \text{ for  }0 < \Theta_1 <1 \,,
\label{WW}
\end{equation}
where $\Theta_1$ is among the set of algorithmic parameters $\Theta$ mentioned above.  This choice of $W$ is symmetric and doubly stochastic. 
Other  choices are given in \cite{boyd2004fastest,yuan2013convergence,tsitsiklis1984problems,chen2012fast,EXTRA34}.
Under our choice of $W$, there are many examples of communication-optimizable algorithms, where the convergence rate depends only on $\tilde{\kappa}(L_{G,w})$.    The convergence rate of EXTRA \cite{shi2015extra,shi2015proximal}  can be explicitly represented, after some straight-forward simplifications by letting $W = W_{L_{G,w}}$ \eqref{WW}:
\begin{equation*}
\frac{ \|\textbf{z}^{k} - \textbf{z}^* \| }{\|\textbf{z}^{0} - \textbf{z}^* \|} \leq \Big( 1+ \min \big \{ \tfrac{1}{p \,\tilde{\kappa}(L_{G,w}) + q }   , \tfrac{1}{r \tilde{\kappa}(L_{G,w})} \big\}  \Big)^{-\tfrac{k}{2}}\,,
\end{equation*}
for some norm $\| \cdot \|$ (where $\textbf{z}$ is related to $\textbf{x}$), and decentralized ADMM \cite{ghadimi2014admm,shi2014linear,mota2013d}:
\begin{equation*}
\frac{ \|\textbf{z}^{k} - \textbf{z}^* \| }{\|\textbf{z}^{0} - \textbf{z}^* \|} \leq \Big( 1+ \sqrt{ p \, \tilde{\kappa}(L_{G,w})^{-2}   + q } - \sqrt{q} \Big)^{-\tfrac{k}{2} }\,,
\end{equation*}
where $p,q,r$ are some coefficients depending only on $\Theta, \mu_f$.
Apparently, while fixing other parameters, a smaller $\tilde{\kappa} (L_{G,w})$ tends to give faster convergence.   We  see the same in extension of EXTRA, PG-EXTRA \cite{shi2015proximal}, for nonsmooth functions.

In order to consider the communication cost for an optimization
procedure, we need the following two quantities.  The first one is the amount of communication per
iteration $U$: 
\begin{equation}
 U (\mu, W) :=\textstyle  \sum_{e \in E } \mu_{e} |W_{e}|_0 ,
\end{equation}
where $\mu = (\mu_{e})$ the communication cost on each edge $e$ and $|\cdot |_0$ as the 1/0 function that indicates if $e$ is used for communication or not.
When $W$ is related to $L_{G,w}$ via \eqref{WW}, we also write $ U
(\mu, L_{G,w} ) := U (\mu, W_{L_{G,w}} )  $ to represent this quantity.  
Next, we introduce the second quantity, the total communication cost for an optimization procedure.  
For this purpose, let us recall that $\Gamma_f$ represents a set of constants depending only on the structure of the target function $f$, and $\Theta$ be a set of parameters not depending on either $f$ of the graph $G$. 
We let $K$ be the number of iterations to achieve the accuracy $\|\textbf{z}^{k} - \textbf{z}^* \|\le \epsilon$,
\begin{equation*}
K=K(\Theta, \Gamma_f,  \tilde{\kappa}(L_{G,w}) , \epsilon),
\end{equation*}
as $K$ clearly depends on $\Gamma_f$ , $\Theta$, $\tilde{\kappa}(L_{G,w}) $, and $\epsilon$. By the definition of $R$ in \eqref{R}, we have
\begin{equation*}
R\big( \Theta, \Gamma_f, \tilde{\kappa}(L_{G,w}), K ( \Theta, \Gamma_f,   \tilde{\kappa}(L_{G,w}) , \epsilon )   \big) = \epsilon \,. 
\end{equation*}
\begin{theorem}
The total communication cost for $ ( \textbf{z}^{i} )_{i=0}^{k}$ generated by a communication-optimizable algorithm to satisfy $\|\textbf{z}^{k} - \textbf{z}^* \| \leq \varepsilon $, denoted as $J( \Theta, \Gamma_f, W_{L_{G,w}} , \epsilon )$, satisfies:
\begin{equation}
J( \Theta, \Gamma_f, W_{L_{G,w}} , \epsilon ) \leq
K \big( \Theta, \Gamma_f,   \tilde{\kappa}(L_{G,w}\big) , \epsilon) \cdot
 U\big (\mu, L_{G,w} \big)  \,.
\label{realcost}
\end{equation}
\end{theorem}
The cost consists of two parts, the first part $K$ decreasing with the connectivity of $G$ while the other part $U$ growing with it.

In general, how to optimize $J$ depends on each algorithm (which gives the function $K$) and the weights $\{ \mu_{e}\}$.
Given an additional assumption that $\{\mu_{e}\}$ are uniformly, we can assume
\begin{equation}
 U(\mu, L_{G,w}) \approx F(d^{\text{ave}}) \,,
\end{equation}
where $d^{\text{ave}} = \sum_i \deg(i) / \sum_i 1$ is the average degree, and $F$ a function that increases monotonically w.r.t. $d^{\text{ave}}$.
Therefore, a reasonable approximation of the communication efficiency optimization problem becomes
\begin{equation}
\min_{G} J( \Theta, \Gamma_f, W_{L_{G,w}} , \epsilon ) \\ \approx \min_{d^{\text{ave}}}\left\{    H (d^{\text{ave}})  F (d^{\text{ave}}) \right\}
\label{happy}
\end{equation}
where $H(d^{\text{ave}})$ reads
\begin{equation}
H(d^{\text{ave}}) :=
 \min_{ ( G_{d^{\text{ave}}},w )}
K \big( \Theta, \Gamma_f,   \tilde{\kappa}(L_{G_{d^{\text{ave}}},w}  ) , \epsilon\big),
\label{happy2}
\end{equation}
where $G_{d^{\text{ave}}}$ denotes an unknown graph whose nodes have the same degree.
Optimizing $d^{\text{ave}}$ in \eqref{happy} after knowing $H (d^{\text{ave}}) $ and $F (d^{\text{ave}}) $ can only be done on a problem-by-problem and algorithm-by-algorithm basis. However, the minimum arguments in the expression \eqref{happy2} in the definition of $H (d^{\text{ave}}) $ are always graphs such that $ \tilde{\kappa}(L_{G_{d^{\text{ave}}},w})$ is minimized.  In light of this, we propose to promote communication efficiency by optimizing $ \tilde{\kappa}(L_{G_{d^{\text{ave}}},w} )$ in any case (i.e. whether we know the actual expression of $K$ or not, or even in the case when $\{\mu_e\}$ is not so uniform.)

\section{Graph optimization (of $\tilde{\kappa}(L_{G})$ or $\tilde{\kappa}(L_{G,w})$)}

\subsection{Exact optimizer with specific nodes and degree $(N,d)$: Ramanujan graphs}

For a general $d$-regular graph, it is known from Alon Boppana Theorem in \cite{alon1997edge,alon2001semi,dodziuk1984difference} that
\begin{equation*}
 \lambda_1 (A) \geq 2 \sqrt{d - 1 } - \frac{ 2 \sqrt{d - 1 } -1}{ \lfloor D/2\rfloor}
\end{equation*}
where $D$ is the diameter of $G$.
\begin{definition} \cite{murty2003ramanujan,newman2000laplacian,chung1997spectral,lubotzky1988ramanujan}
A $d$-regular graph $G$ is a Ramanujan graphs if it satisfies
\begin{equation}
 \lambda_1 (A_G) \leq 2 \sqrt{d - 1 },
\end{equation}
where $A_G$ is the adjacency matrix of $G$.
\end{definition}

In fact, a Ramanujan graph serves as an asymptotic minimizer of $\tilde{\kappa}(L_{G_d}) $ for a $d$-regular graph $G_d$. For a Ramanujan graph, the reduced condition number of the Laplacian satisfies:
\begin{equation*}
\tilde{\kappa}(L_{G_d})  \leq \frac{d + 2 \sqrt{d - 1 } }{ d - 2 \sqrt{d - 1 } } \,.
\label{eqLd}
\end{equation*}


If we use a Ramanujan graph that minimizes $\tilde{\kappa}(L_{G_d})$, we can then find an approximate minimizer in \eqref{happy2} and thus \eqref{happy} for the total communication cost by further finding $d$ (or $d^{\text{ave}}$).

Explicit constructions of Ramanujan graphs are known only for some special $(N,d)$, as labelled Cayley graphs and the Lubotzky-Phillips-Sarnak construction \cite{lubotzky1986explicit}, with some computation given in \cite{kpsgraph,randy2010} or construction of bipartite Ramanujan graphs
\cite{marcus2013interlacing}.
Interested readers may consult \cite{newman2000laplacian,murty2003ramanujan,hoory2006expander,chung1997spectral,lubotzky1986explicit,lubotzky1988ramanujan} for more references.

\subsection{Optimizer for large $N$: Random $d$-regular graphs}

In some practical situations, however, we encounter pairs $(N,d)$ where an explicit Ramanujan graph is still unknown. When $N$ is very large, we propose a random $d$-regular graph as a solution of optimal communication efficiency, with randomness performed as follows \cite{friedman2008proof} (where $d$ shall be even and $N$ can be an arbitrary integer that is greater than $2$):
choosing independently $d/2$ permutations $\pi_{j}, j = 1,..,d/2 $ of the numbers from $1$ to $n$, with each permutation equally
likely, a graph $G = (V,E)$ with vertices $V = {1,...,n}$ is given by
\begin{equation*}
E = \{ (i, \pi_{j}(i) ), i = 1,...n, j = 1,..., d/2 \}\,.
\end{equation*}
The validity of using random graphs as an approximate optimizer is ensured by the Friedman Theorem as below:
\begin{theorem} \cite{friedman2008proof}
For every $\epsilon>0$,
\begin{equation}
\mathbb{P} (\lambda_1(A) \leq 2 \sqrt{d-1} + \epsilon)
= 1- o_N(1)
\end{equation}
where $G$ is a random $(N,d)$-graph.\footnote{The notion $g(N)=o_N(f(N))$ means $\lim_{N\rightarrow\infty}g(N)/f(N)=0$.}
\end{theorem}
Roughly speaking, this theorem says that for a very large $N$, almost all random $(N,d)$ graphs are Ramanujan, i.e. they are nicely connected. Therefore, it is just fine to adopt a random $(N,d)$ graph when $N$ is too large for a Ramanujan graph or sparsifier to be explicitly computed.

It is possible that some regular graphs are not connected. To address this issue, we first grow a random spanning tree of the underlying graph to ensure connectivity. A random generator for regular graphs is  GGen2 \cite{GGen}. 


\subsection{Approximation for small $N$: $2$-Ramanujan Sparsifier} \label{sec:sparsifier}
For some small $N$, an explicit $(N,d)$-Ramanujan graph may still be unavailable. We hope to construct an approximate Ramanujan graph, i.e. expanders which are sub-optimal but work well in practice.  An example is given by the $2$-Ramanujan sparsifier, which is the subgraph $H$ as follows:

\begin{theorem} \cite{batson2012twice}
\label{imp_sparse}
For every $d >1$, every undirected weighted graph $G = (V,E,w)$ on $N$ vertices contains a weighted subgraph
$H = (V,F, w')$ with at most $d(N-1)$ edges (i.e. an average at most $2d$ graph) that satisfies
\begin{equation}
L_{G,w} \preccurlyeq L_{H,w'} \preccurlyeq \frac{d+1 + 2 \sqrt{d}}{d+1 - 2 \sqrt{d}}  \, L_{G,w}
\label{hahahaha}
\end{equation}
where the Loewner partial order relation $A \preccurlyeq B$ indicates that $B-A$ is a positive semidefinite matrix.
\end{theorem}
This theorem allows us to construct a sparsifier from an original graph $G$.  Actually, the proof of the theorem provides us with an explicit algorithm for such a construction \cite{batson2012twice}.

\section{Numerical Experiments} \label{sec:num}

In this subsection, we illustrate the communication efficiency optimization achieved by expander graphs.

Since the focus of our work is not to investigate various methods to produce expander graphs, we did not test \textbf{Algorithm I}.
Rather, our focus is to illustrate that a smaller \emph{reduced condition number} improves communication efficiency and that expander graphs are good choices to reduce communication in decentralized optimization.  Therefore we  compare the difference in communication efficiency using existing graphs with different reduced condition numbers.
We compare the convergence speeds on a family of (possibly non-convex non-smooth) problems which are approximate formulations of finding the sup-norm of a vector $l = (l_1,...,l_M)$ (up to sign).

\noindent \textbf{Example 1}.  We choose $M = P = 1092$.  Our problem is:
\begin{equation*}
\min_{x \in \mathbb{R}^P } f(x) = \tfrac{1}{M}\textstyle
 \sum_{i= 1}^M  f_{i}(x) ,
\end{equation*}
where
\begin{equation*}
f_i(x) = l_i x_i /  ( \| x \|_1 + \varepsilon )
\end{equation*}
for a small $\varepsilon$ chosen as  $\varepsilon = 1 \times 10^{-5}.$
One very important remark is that this optimization problem is non-convex and non-smooth, and the $\varepsilon$-optimizer occupies a tabular neighbourhood of a whole ray $\{ -\lambda e_{n} : \lambda > 0 \}$, where $l_n = \| l \|_{\infty}$ and $e_n$ is the $n$th coordinate vector $(0,0,...,1,....,0)$.

\begin{figure}
\begin{center}
     \vskip -0.3truecm
           \scalebox{0.35}{\includegraphics{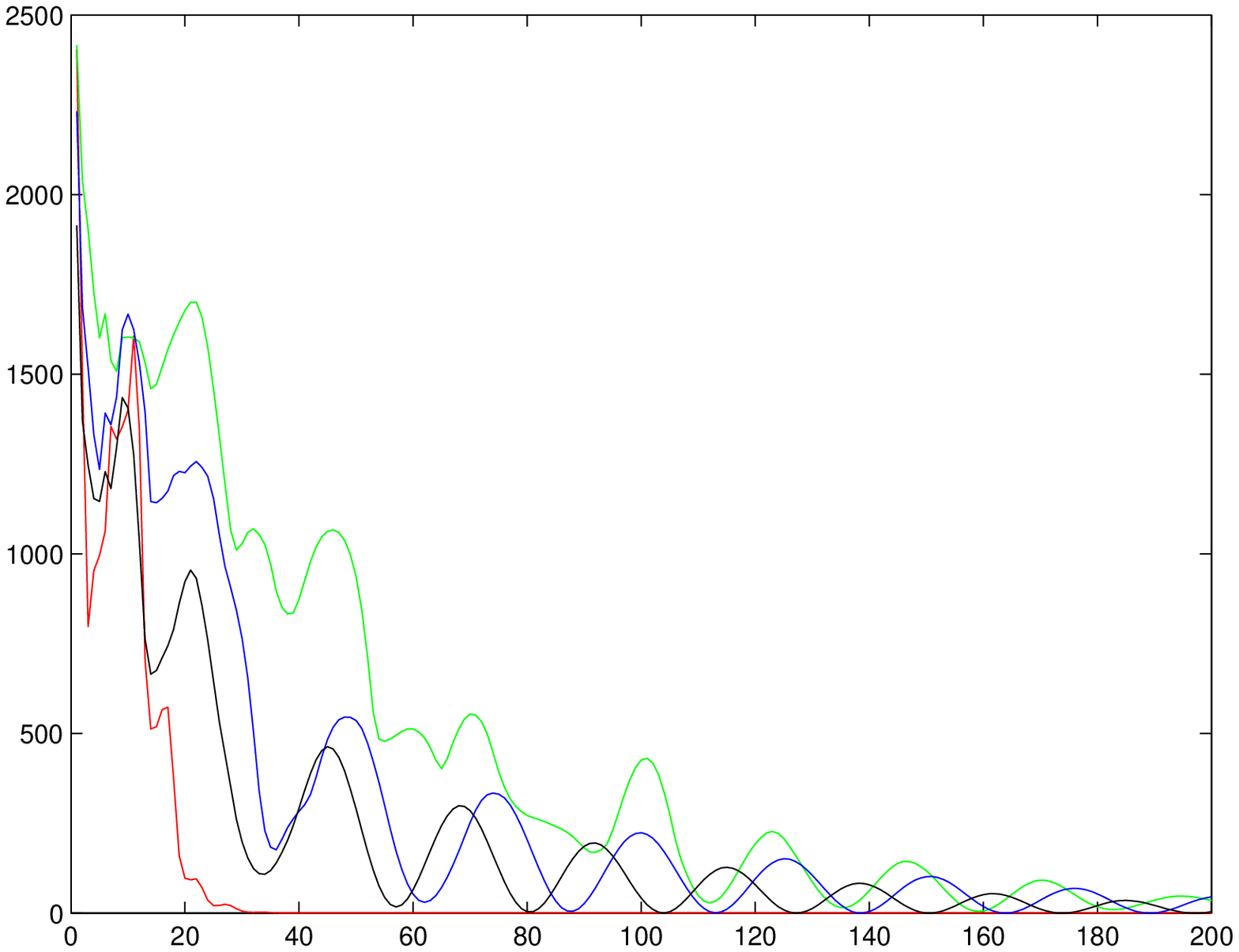}}
    \vskip -0.5truecm
    \scalebox{0.35}{\includegraphics{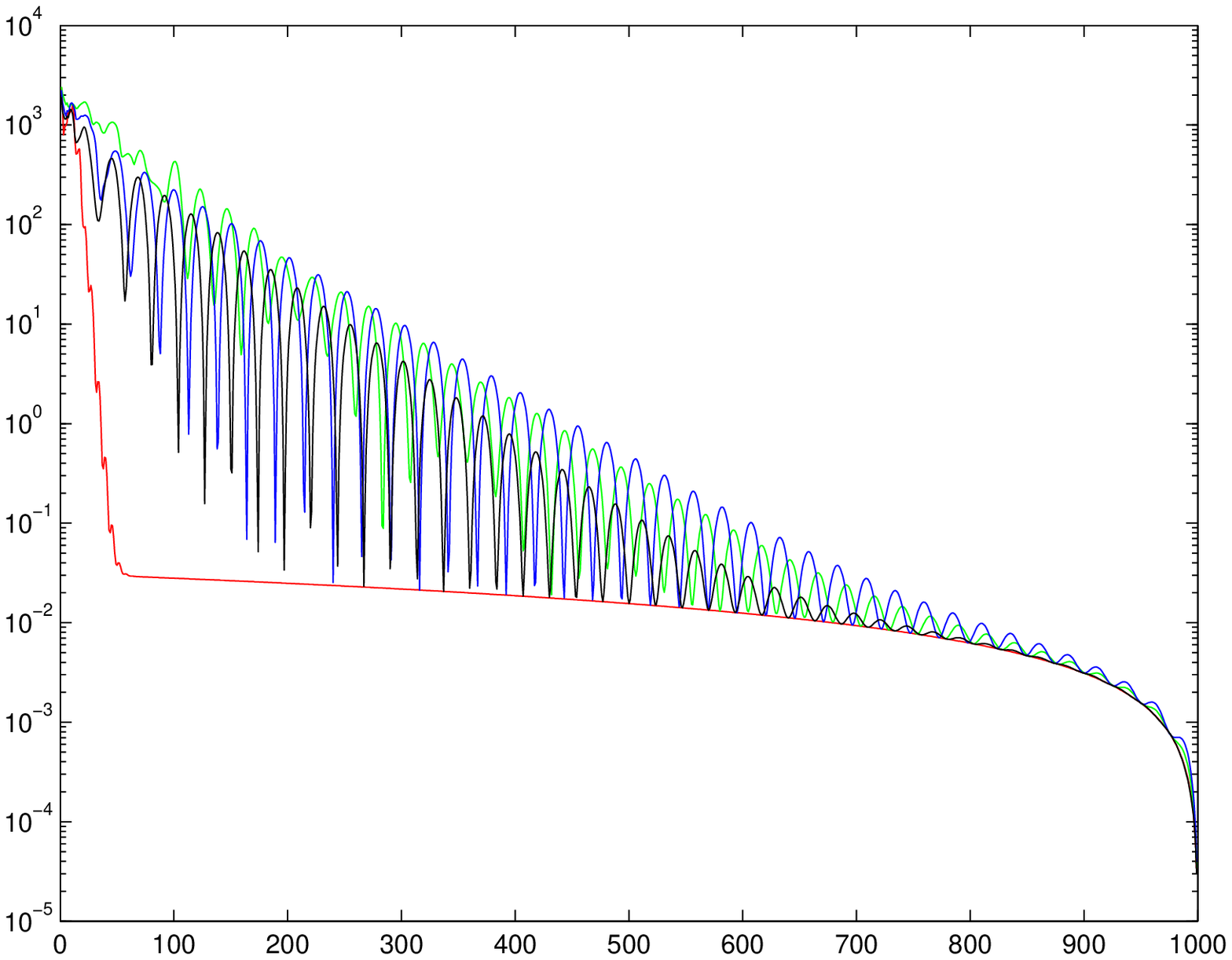}}
    \vskip -0.3truecm
     \caption{\small Convergence rate of EXTRA in \textbf{Example 1} over Ramanujan $30$-graph $\LPS(29,13)$ (red, the best), $30$-reg graph (blue), $60$-reg graph (green) and $120$-reg graph (black).}\label{EXTRA}
     \end{center}
\end{figure}

We solve the problem using the EXTRA algorithm \cite{shi2015extra}:
\begin{equation*}
\begin{cases}
\textbf{x}^{1} = W \textbf{x}^{0} - \alpha  \nabla \textbf{f} (\textbf{x}^{0})  \,, \\
\textbf{x}^{k+2} = (W+I ) \textbf{x}^{k+1} - \frac{W+I }{2} \textbf{x}^{k} - \alpha [\nabla \textbf{f} (\textbf{x}^{k+1}) - \nabla \textbf{f} (\textbf{x}^{k})  ] \,.
\end{cases}
\end{equation*}
To compare speeds, we plot the following sequence:
\begin{equation*}
  \delta_k := \tfrac{1}{M} \textstyle\sum_{i= 1}^M   F_{i}( \textbf{x}^k ) - \min_{x} F(x) \,.
\end{equation*}
Figure \ref{EXTRA} shows the sequence $ \delta_k$ under $4$ different graph topologies:
Ramanujan $30$-graph $\LPS(29,13)$, a regular-$30$, a regular-$60$ graph, and a regular-$120$ graph. 
Our regular $d$ graphs (other than $LPS(29,13)$) are generated by joining the $i$th node to the $\left \{ i+ \lfloor \frac{N}{d} \rfloor k  \mod  N \, : \, 0 \leq k < \frac{d}{2} \right\}$th nodes for all $ 0 \leq i \leq N-1 $ (in here we label the nodes from $0$ to $N-1$).
The mixing matrices $W$ that we use are generated as $W_{L_{G,w}}$ as described in \eqref{WW} with $w = (w_e )$ where $w_e = 1$ on all edges $e \in E$.
We clearly see that the convergence rate of the algorithm under Ramanujan $30$-graph $\LPS(29,13)$ is even better than that of a regular-$120$ graph in the first $40$ iterations.   Afterward, the curve of the Ramanujan $30$-graph becomes flat, probably because it already arrived a small tabular neighborhood of the ray of $\varepsilon$-optimizer. Table \ref{EXTRA_table} shows the number of iteration $k_0$ as well as its total communication complexity to achieve $$ \text{stopping rule:}\quad |\delta_{k_0}  -  \delta_{k_0-1}| < 1\times 10^{-3}.$$ 
We see that an expander graph reduces communication complexity.  One can observe the rapid drop of the red curve because of its efficiency to locate the active set.

\begin{table}
\begin{tabular}{ c| r | r | r}
\hline
 & $\tilde{\kappa}(L_G)$ & No. of Itr. &             Total comm. complexity \\
\hline
$\LPS(29,13)$  & 1.9538 &          52  &1703520 \\
\hline
Regular $30$ graph &  30.5375  &     444  & 14545440 \\
\hline
Regular $60$ graph &   32.1103   &     494  & 32366880 \\
\hline
Regular $120$ graph &   27.1499 &      454  &  59492160 \\
\hline
\end{tabular}
 \caption{\small Communication complexity in \textbf{Example 1} for Ramanujan $30$-graph $\LPS(29,13)$, $30$-reg, $60$-reg, and$120$-reg graph.}\label{EXTRA_table}
\end{table}

\begin{figure}
\begin{center}
     \vskip -0.3truecm
           \scalebox{0.35}{\includegraphics{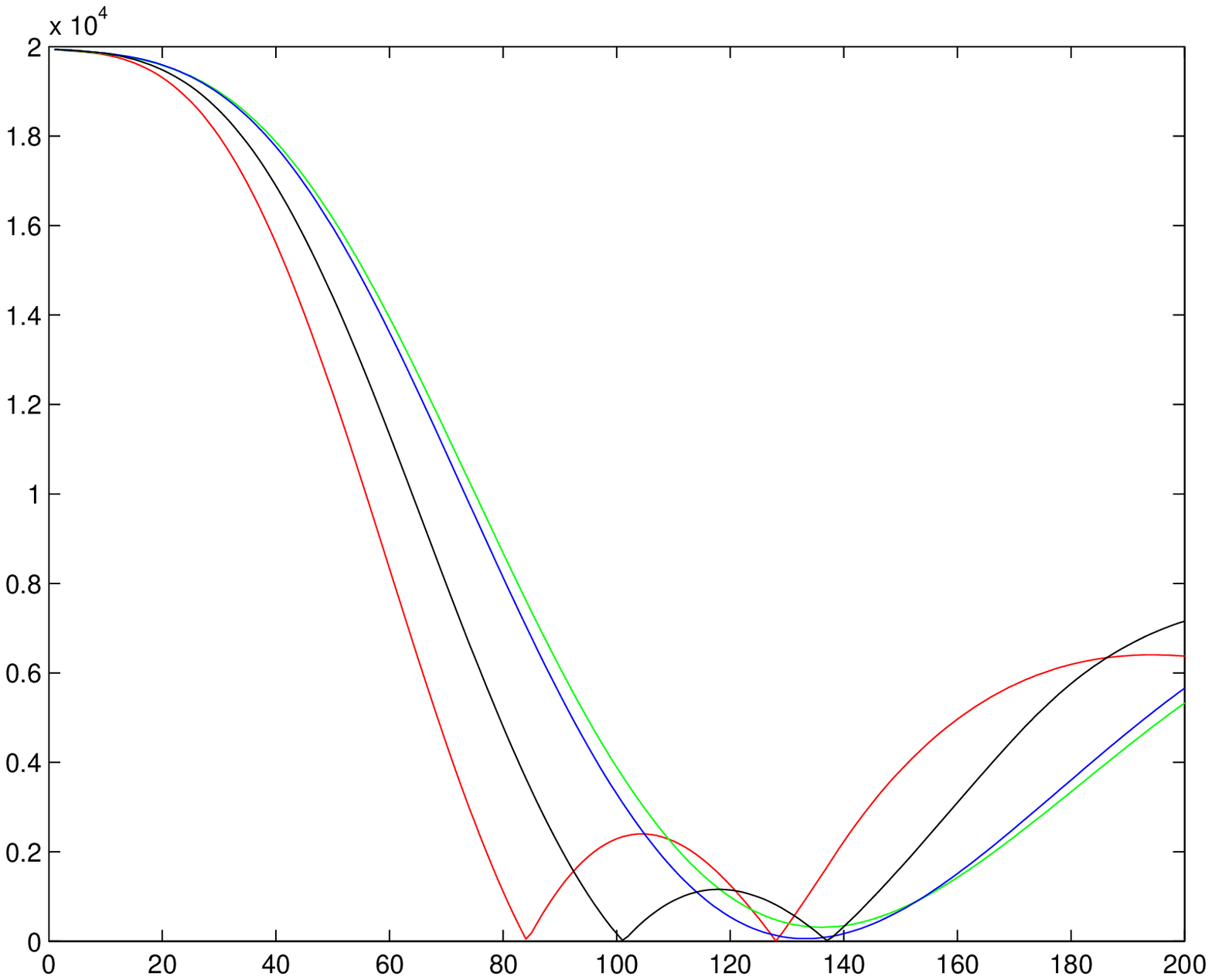}}
    \vskip -0.5truecm
    \scalebox{0.35}{\includegraphics{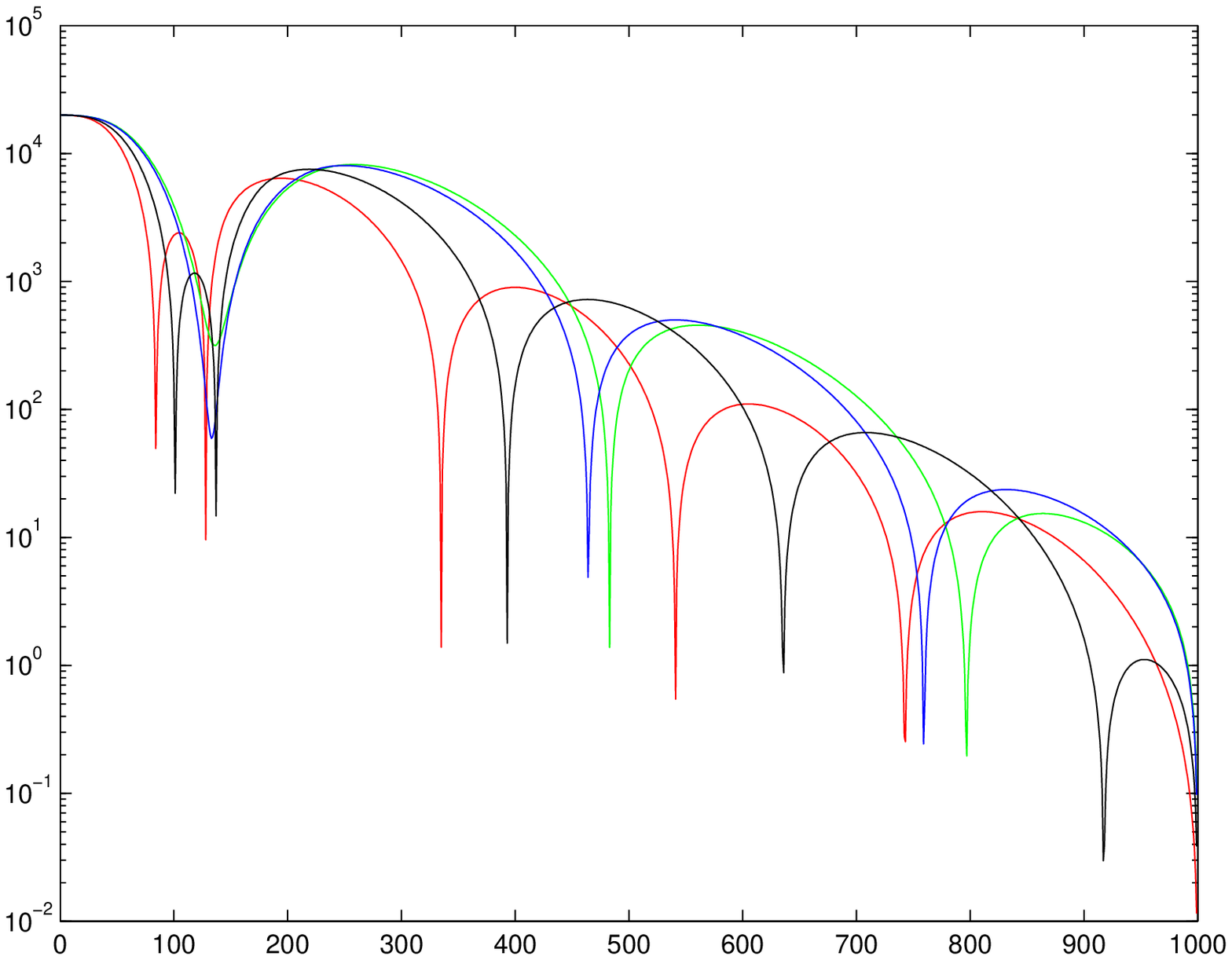}}
    \vskip -0.3truecm
        \caption{\small Convergence rate of PG-EXTRA in \textbf{Example 2} over Ramanujan $30$-graph $\LPS(29,13)$ (in red), $30$-reg (in blue), $60$-reg (in green), and $120$-reg graph (in black).}\label{PGEXTRA2}
     \end{center}
\end{figure}

\noindent \textbf{Example 2}.  In this case we again choose $M = 1092$, but $P = 1 $ (i.e. $1$ dimensional problem).  Our problem is
with a similar form as in \textbf{Example 1}, but with
\begin{equation*}
f_i(x) = \chi_{a \geq \sqrt { | l_i | } }(x ) + |x |^2.
\end{equation*}
This is a convex non-smooth problem. We apply  PG-EXTRA \cite{shi2015proximal} to solve this problem.
To proceed, we split $f_i$ as
\begin{equation*}
\begin{cases}
s_i(x) := |x |^2  \\
r_i(x) :=  \chi_{a \geq \sqrt { | l_i | } }(x ) ,
\end{cases}
\end{equation*}
and apply the PG-EXTRA iteration:
\begin{equation*}
\begin{cases}
\textbf{x}^{\frac{1}{2}} = W \textbf{x}^{0} - \alpha  \nabla \textbf{f} (\textbf{x}^{k+1})  \,, \\
\textbf{x}^{1} = \argmin \textbf{r}(\textbf{x}) + \frac{1}{2 \alpha} \| \textbf{x} - \textbf{x}^{\frac{1}{2}}  \|_2^2 \,,\\
\textbf{x}^{k+\frac{3}{2}} = W \textbf{x}^{k+1} + \textbf{x}^{k+\frac{1}{2}} - \frac{W+I }{2} \textbf{x}^{k} - \alpha [\nabla \textbf{s} (\textbf{x}^{k+1}) - \nabla \textbf{s} (\textbf{x}^{k})  ] \,, \\
\textbf{x}^{k+2} = \argmin \textbf{r}(\textbf{x}) + \frac{1}{2 \alpha} \| \textbf{x} - \textbf{x}^{k + \frac{3}{2}}  \|_2^2 \,.\\
\end{cases}
\end{equation*}
To compare speeds, we plot the following sequence:
\begin{equation*}
  \delta_k := \big| \tfrac{1}{M} \textstyle\sum_{i= 1}^M  s_{i}( \textbf{x}_i^k ) - \min_{x} F(x)  \big| \,.
\end{equation*}
and exclude the $r_i(x) $ part.
The mixing matrices that we use are the same as those described in \textbf{Example 1}.
As shown in Figure \ref{PGEXTRA2} the convergence rate of the algorithm under the Ramanujan $30$-graph $\LPS(29,13)$ is still the best though is  less outstanding.
Table \ref{PGEXTRA2_table} shows the numbers of iteration $k_0$ as well as the total communication complexities to achieve $  | \delta_{k_0}  -  \delta_{k_0-1} | < 1\times 10^{-3}$. The expander graph has the clear advantage.
\begin{table}
\begin{tabular}{ c| r | r | r}
\hline
 & $\tilde{\kappa}(L_G)$ & No. of Itr. &             Total comm. complexity \\
\hline
$\LPS(29,13)$  & 1.9538 &          811 &  26568360 \\
\hline
Regular $30$ graph &  30.5375  &     1001  &  32792760 \\
\hline
Regular $60$ graph &  32.1103    &     832  & 54512640 \\
\hline
Regular $120$ graph &  27.1499 &     953 &   124881120 \\
\hline
\end{tabular}
 \caption{\small Communication complexity in \textbf{Example 2} for Ramanujan $30$-graph $\LPS(29,13)$, $30$-reg, $60$-reg  and$120$-reg graphs.}\label{PGEXTRA2_table}
\end{table}

\bibliographystyle{IEEEbib}
\bibliography{mybibfile}

\end{document}